\documentclass[11pt]{amsart}
\usepackage{amssymb,booktabs,graphicx}
\newcommand{\R}{\mathbb{R}}
\newcommand{\Z}{\mathbb{Z}}
\newcommand{\Q}{\mathbb{Q}}
\newcommand{\norm}[1]{\lVert #1\rVert}
\newcommand{\hash}[2]{\texttt{#1}\newline\texttt{#2}}
\newtheorem{theorem}{Theorem}
\newtheorem{lemma}[theorem]{Lemma}
\title[Five colours for six regular polygon norms]{The chromatic number of the Minkowski plane\\ with a regular 14-, 16-, 18-, 20-, 22- or 24-gon\\ as unit ball is at least five}
\author{King Tat Wong}
\email{info@amberlogydesignhouse.com}
\date{}
\begin{document}
\begin{abstract}
Let $P_{2n}$ be the regular $2n$-gon centred at the origin.
We give finite unit-distance graphs, none of them $4$-colourable, in the norms with unit balls $P_{14}$, $P_{16}$, $P_{18}$, $P_{20}$, $P_{22}$ and $P_{24}$.
Thus $\chi(\mathbb{R}^2,P_{2n})\ge 5$ in these six cases.
With Geh\'er's upper bound this yields $\chi\in\{5,6\}$ for $2n\in\{14,16,18,20,22\}$; for $2n=24$ the upper bound we use is Chilakamarri's $\chi\le 7$.
A second graph in the $22$-gon norm, with $6517$ vertices and $42394$ edges, is certified by the same chain and is not the graph used for the theorem.
Non-$4$-colourability is certified by DRAT proofs checked with \texttt{drat-trim} and by LRAT proofs checked with \texttt{lrat-check}.
\texttt{cake\_lpr} accepts every LRAT proof of the six graphs giving the bound above, and of the second $22$-gon graph.
Edges are checked by exact arithmetic, and also in Lean~4.
The $16$-gon interval check uses a separation lemma proved in full below.
The graphs in the theorem range from $225$ to $5577$ vertices.
The certificates for $2n\in\{14,18,20,22,24\}$ are at doi:10.5281/zenodo.24458454, together with the Lean sources for $2n\in\{14,16,18,20,22,24\}$ and for the second $22$-gon graph.
The $16$-gon certificate files are at doi:10.5281/zenodo.24458031.
\end{abstract}
\maketitle

\section{Introduction}
Exoo, Fisher and Ismailescu proved $\chi(\mathbb{R}^2,P_{2n})\ge 5$ for $2n=8,10,12$ and asked whether the same holds for every regular $2n$-gon with $n\ge 4$ \cite[Problem~5.1]{EFI}.
Geh\'er proved $\chi\le 6$ for every regular $2n$-gon with $2n\le 22$ \cite{Geher}.
For larger $2n$, Chilakamarri's bound $\chi\le 7$ is the one we use \cite{Chil}.

\begin{theorem}\label{thm:main}
$\chi(\mathbb{R}^2,P_{2n})\ge 5$ for $2n\in\{14,16,18,20,22,24\}$.
The graphs are: $225$ vertices and $1176$ edges for $2n=14$; $960$ vertices and $5248$ edges for $2n=16$; $612$ vertices and $3042$ edges for $2n=18$; $1861$ vertices and $11615$ edges for $2n=20$; $5577$ vertices and $40603$ edges for $2n=22$; $1296$ vertices and $8616$ edges for $2n=24$.
Consequently $\chi(\mathbb{R}^2,P_{2n})\in\{5,6\}$ for $2n\in\{14,16,18,20,22\}$, and $5\le\chi(\mathbb{R}^2,P_{24})\le 7$.
\end{theorem}

Section~\ref{sec:norm} fixes the norm, the colouring formula and the separation bounds.
Section~\ref{sec:chain} records how those pieces give the theorem.
The graphs are in Section~\ref{sec:graphs}, and the certificates are in Section~\ref{sec:data}.

\section{Norm, colouring, and a separation bound}\label{sec:norm}
Let $\zeta=e^{i\pi/n}$ and $P_{2n}=\operatorname{conv}\{\zeta^k:0\le k<2n\}$.
Write $L_k(z)=\operatorname{Re}(z\bar\zeta^k)+\operatorname{Re}(z\bar\zeta^{k+1})$ and $h=1+\cos(\pi/n)$.

\begin{lemma}[Unit distance]\label{lem:unit}
$\|z\|_{P_{2n}}=\max_k L_k(z)/h$.
\end{lemma}
\begin{proof}
$L_k/h$ equals $1$ on the side $[\zeta^k,\zeta^{k+1}]$, and $P_{2n}$ is the intersection of the half-planes $L_k\le h$.
\end{proof}

\begin{lemma}[Symmetry breaking]\label{lem:tri}
If $a,b,c$ are pairwise adjacent, then a graph is $4$-colourable if and only if it has a proper $4$-colouring in which $a,b,c$ receive colours $0,1,2$.
\end{lemma}
\begin{proof}
A proper colouring gives the triangle three distinct colours. Permute the colours.
\end{proof}

For a graph on vertices $0,\dots,N-1$ the colouring formula uses variables $x_{v,c}$ ($c\in\{0,1,2,3\}$), numbered $4v+c+1$, and the clauses $\bigvee_c x_{v,c}$, $\neg x_{v,a}\vee\neg x_{v,b}$ for $a<b$, and $\neg x_{u,c}\vee\neg x_{v,c}$ for each edge $\{u,v\}$ and each colour $c$.
Satisfying assignments are exactly the proper $4$-colourings.
In every triangle certificate below, one triangle is fixed to colours $0,1,2$ by three unit clauses.
Lemma~\ref{lem:tri} says that this does not change satisfiability.

\begin{lemma}[Separation]\label{lem:sep}
Let $d=c/\mathrm{sc}$ with $c=\sum_j c_j\zeta^j\in\mathbb{Z}[\zeta]$, $\mathrm{sc}$ a positive integer, $m=[K:\mathbb{Q}]$ where $K=\mathbb{Q}(\zeta)\cap\mathbb{R}$, and $B=4\sum_j|c_j|+4\,\mathrm{sc}$.
If $n\ge 3$ and $\|d\|_{P_{2n}}\ne 1$, then $\bigl|\|d\|_{P_{2n}}-1\bigr|>1/(4\,\mathrm{sc}\,B^{m-1})$.
\end{lemma}
\begin{proof}
Write $\|{\cdot}\|=\|{\cdot}\|_{P_{2n}}$ and let $k$ attain the maximum, so $\|d\|=L_k(d)/h$.
Set $\beta=2\,\mathrm{sc}\,(L_k(d)-h)$; then $\|d\|-1=\beta/(2\,\mathrm{sc}\,h)$.
Since $n\ge 3$, one has $0<\pi/n<\pi/2$, so $0<\cos(\pi/n)<1$ and $1<h<2$.
Thus if $|\beta|\ge B^{-(m-1)}$, then
\[
  \bigl|\|d\|-1\bigr|=\frac{|\beta|}{2\,\mathrm{sc}\,h}>\frac{|\beta|}{4\,\mathrm{sc}}\ge\frac{1}{4\,\mathrm{sc}\,B^{m-1}}.
\]
It remains to check the lower bound on $|\beta|$, and that $\beta\ne 0$ when $\|d\|\ne 1$.
The second is immediate from $h>0$.

Now $2\operatorname{Re}(c\bar\zeta^k)=c\bar\zeta^k+\bar c\,\zeta^k\in\mathbb{Z}[\zeta]$, and likewise for $k+1$.
Also $2h=2+\zeta+\zeta^{-1}$, so $2\,\mathrm{sc}\,h\in\mathbb{Z}[\zeta]$.
Hence $\beta\in\mathbb{Z}[\zeta]$.
Complex conjugation fixes $\beta$, so $\beta$ lies in $K$ and is an algebraic integer of $K$.

Let $\sigma$ be an automorphism of $\mathbb{Q}(\zeta)$.
Then $|\sigma(\zeta^j)|=1$, so $|\sigma(c\bar\zeta^k)|\le\sum_j|c_j|$ and $|\sigma(\bar c\,\zeta^k)|\le\sum_j|c_j|$.
Therefore $|\sigma(2\operatorname{Re}(c\bar\zeta^k))|\le 2\sum_j|c_j|$, and the same bound holds for the term with $k+1$.
Also $|\sigma(2\,\mathrm{sc}\,h)|=|\sigma(\mathrm{sc}(2+\zeta+\zeta^{-1}))|\le\mathrm{sc}(2+1+1)=4\,\mathrm{sc}$.
Adding these estimates gives $|\sigma(\beta)|\le B$.
If $\beta\ne 0$, the field norm $N_{K/\mathbb{Q}}(\beta)$ is a non-zero rational integer, so $|N|\ge 1$.
The other $m-1$ conjugates each have absolute value at most $B$, hence $|\beta|\ge B^{-(m-1)}$.
\end{proof}

Three programs check the edges, and they share no code.
Method~(1) compares support functionals by exact arithmetic in $\Z[\alpha]$.
Method~(2) writes the difference in sector coordinates and decides the norm by exact algebraic signs.
Method~(3) evaluates the norm in decimal interval arithmetic.
For $2n\in\{14,18,20,22,24\}$, Lemma~\ref{lem:sep} is used only by method~(3).
Methods~(1) and~(2), and the Lean proofs, decide equality exactly.
For the $16$-gon, $\alpha=2\cos(\pi/8)$ satisfies $\alpha^2=2+\sqrt2$ and $\sin(\pi/8)=(\sqrt2-1)\alpha/2$.
Write $\det(a,b)=\operatorname{Im}(\bar a b)$, $V_k=2\alpha\zeta^k$, and $Q_{16}=\operatorname{conv}\{V_k:0\le k\le 15\}=2\alpha P_{16}$ for the scaled unit ball in which the $16$-gon file coordinates are stored.
Then
\[
  \norm{d}_{Q_{16}}=\norm{d/(2\alpha)}_{P_{16}}=\max_{0\le k\le 15} f_k(d),\qquad
  f_k(d)=\frac{\det(d,\,V_{k+1}-V_k)}{\det(V_k,V_{k+1})}.
\]
The functional $f_k$ equals $1$ on the line through $V_k$ and $V_{k+1}$.
Methods~(1) and~(2) decide whether $f_k(d)=1$ by exact arithmetic.
Method~(3) uses interval arithmetic with the following sharper bound.

\begin{lemma}[Separation for $Q_{16}=2\alpha P_{16}$]\label{lem:gap16}
Let $d=c/2$ with $c=\sum_{j=0}^7c_j\zeta^j\in\Z[\zeta]$, write $\norm{\cdot}=\norm{\cdot}_{Q_{16}}$, and put
$B=4\sum_j|c_j|+6$. If $\norm{d}\ne1$, then
$\bigl|\norm{d}-1\bigr|>1/(6B^3)$.
\end{lemma}

\begin{proof}
Write $\det(a,b)=\operatorname{Im}(\bar a b)$. With $V_k=2\alpha\zeta^k$ and $\alpha\in\R$,
\[
  \overline{V_k}\,V_{k+1}=(2\alpha)^2\zeta^{-k}\zeta^{k+1}=4\alpha^2\zeta,
\]
so $\det(V_k,V_{k+1})=4\alpha^2\sin(\pi/8)$. The identities
$\sin(\pi/8)=(\sqrt2-1)\alpha/2$ and $\alpha^2=2+\sqrt2$ give
\begin{align*}
  4\alpha^2\sin(\pi/8)
  &=4\alpha^2\cdot\frac{(\sqrt2-1)\alpha}{2}
  =2\alpha^3(\sqrt2-1)\\
  &=2\alpha\cdot\alpha^2(\sqrt2-1)
  =2\alpha\cdot(2+\sqrt2)(\sqrt2-1).
\end{align*}
Expanding the product,
\[
  (2+\sqrt2)(\sqrt2-1)
  =2\cdot\sqrt2+\sqrt2\cdot\sqrt2-2\cdot 1-\sqrt2\cdot 1
  =2\sqrt2+2-2-\sqrt2
  =\sqrt2,
\]
hence $\det(V_k,V_{k+1})=2\alpha\cdot\sqrt2=2\sqrt2\,\alpha$.

Let $w=\zeta^{k+1}-\zeta^k$ and $z=\bar c\,w$, and let $k$ attain the maximum,
so $\norm{d}=f_k(d)$. Since $d=c/2$ and $V_{k+1}-V_k=2\alpha w$,
\begin{align*}
  \det(d,\,V_{k+1}-V_k)
  &=\operatorname{Im}\bigl(\overline{(c/2)}\cdot 2\alpha w\bigr)
   =\operatorname{Im}(\alpha z)
   =\alpha\operatorname{Im}(z),
\end{align*}
the last step because $\alpha$ is real. Dividing by
$\det(V_k,V_{k+1})=2\sqrt2\,\alpha$ therefore gives
\[
  f_k(d)
  =\frac{\alpha\operatorname{Im}(z)}{2\sqrt2\,\alpha}
  =\frac{\operatorname{Im}(z)}{2\sqrt2}.
\]
From $\operatorname{Im}(z)=(z-\bar z)/(2i)$ and $1/i=-i$,
\[
  \frac{z-\bar z}{2i}=\frac{-i}{2}(z-\bar z),
\]
so $\operatorname{Im}(z)=-i(z-\bar z)/2$ and
\[
  f_k(d)
  =\frac{-i(z-\bar z)/2}{2\sqrt2}
  =\frac{-i(z-\bar z)}{4\sqrt2}.
\]
Hence
\begin{align*}
  \beta
  &:=4\sqrt2\bigl(f_k(d)-1\bigr)
  =4\sqrt2\left(\frac{-i(z-\bar z)}{4\sqrt2}-1\right)\\
  &=-i(z-\bar z)-4\sqrt2.
\end{align*}

Here $i=\zeta^4$. Also $\zeta^2=e^{i\pi/4}=(\sqrt2/2)(1+i)$ and
$\zeta^6=e^{i3\pi/4}=(\sqrt2/2)(-1+i)$, so
\[
  \zeta^2-\zeta^6=\frac{\sqrt2}{2}\bigl((1+i)-(-1+i)\bigr)=\sqrt2.
\]
Since $\bar\zeta=\zeta^{-1}$, one has $\bar c=\sum_j c_j\zeta^{-j}\in\Z[\zeta]$,
hence $z\in\Z[\zeta]$ and $\beta\in\Z[\zeta]$. Conjugation sends $-i$ to $i$ and
$z-\bar z$ to $\bar z-z$, so
\[
  \overline{-i(z-\bar z)}=i(\bar z-z)=-i(z-\bar z),
\]
while $\sqrt2$ is real. Thus conjugation fixes $\beta$, and $\beta$ lies in the
maximal real subfield $K$. As an element of $\Z[\zeta]$ it is an algebraic
integer of $K$.

If $\norm{d}\neq 1$, then $f_k(d)\neq 1$ and $\beta\neq 0$. Every automorphism
$\sigma$ of $\Q(\zeta)$ satisfies $\sigma(\zeta)=\zeta^m$ for some odd $m$, and
the four conjugates of $\beta\in K$ are restrictions of such automorphisms. For
each of them $|\sigma(\zeta^j)|=1$, so $|\sigma(w)|\le 2$ and
$|\sigma(\bar c)|\le\sum_j|c_j|$. Therefore
\[
  |\sigma(z)|=|\sigma(\bar c)\,\sigma(w)|\le 2\sum_j|c_j|.
\]
The complex conjugate is $\bar z=c\,\bar w$, and the same estimate gives
$|\sigma(\bar z)|\le 2\sum_j|c_j|$. Also $|\sigma(i)|=|\sigma(\zeta^4)|=1$.
Since $\sigma$ fixes $\Q$, one has $\sigma(\sqrt2)^2=\sigma(2)=2$, so
$\sigma(\sqrt2)=\pm\sqrt2$ and $|\sigma(\sqrt2)|=\sqrt2$. Hence
\begin{align*}
  |\sigma(\beta)|
  &= \bigl|-\sigma(i)\bigl(\sigma(z)-\sigma(\bar z)\bigr)-4\sigma(\sqrt2)\bigr| \\
  &\le |\sigma(z)|+|\sigma(\bar z)|+4\sqrt2
   \le 4\sum_j|c_j|+4\sqrt2.
\end{align*}
Both $4\sqrt2$ and $6$ are positive, and $(4\sqrt2)^2=32<36=6^2$, so
$4\sqrt2<6$. Therefore
\[
  |\sigma(\beta)|<4\sum_j|c_j|+6=B.
\]
The field norm $N_{K/\Q}(\beta)$ is a non-zero rational integer, so
$|N_{K/\Q}(\beta)|\ge 1$. Each of the other three conjugates has absolute value
strictly less than $B$, and
\[
  |\beta|
  =\frac{|N_{K/\Q}(\beta)|}{\prod_{\sigma\neq\mathrm{id}}|\sigma(\beta)|}
  >\frac{1}{B^3}.
\]
In particular $|\beta|\ge 1/B^3$. Together with $4\sqrt2<6$ this yields
\[
  \bigl|f_k(d)-1\bigr|
  =\frac{|\beta|}{4\sqrt2}
  >\frac{1}{4\sqrt2\,B^3}
  >\frac{1}{6B^3}.
\]
Since $\norm{d}=f_k(d)$, this is the claim.
\end{proof}

\section{From a finite graph to the bound}\label{sec:chain}
A finite unit-distance graph with no proper $4$-colouring forces $\chi\ge 5$.
Each case below is read in the same way.

Every recorded edge has norm exactly $1$.
Methods~(1) and~(2), where they are run, decide this by exact arithmetic, and so does Lean for every graph below.
Method~(3) only checks that the computed norm lies inside the separation gap of $1$.
For $2n\in\{14,16,18,20,24\}$, and for both $22$-gon graphs, all three methods were run, and each accepted every edge.
On the $14$-gon and on each $22$-gon graph, method~(2) also places every vertex in $U+U$.
It was run on a copy of the graph file with the side parameters written after the scale line; the hashed graph files contain no parameter lines.

The colouring formula is the CNF defined above, built from the graph file, with one triangle fixed to colours $0,1,2$.
Lemma~\ref{lem:tri} says that fixing the triangle does not change satisfiability.
CaDiCaL \cite{CaDiCaL} writes a DRAT refutation of that formula.
\texttt{drat-trim} \cite{DRAT} checks the refutation and writes an LRAT proof, which \texttt{lrat-check} accepts.
\texttt{cake\_lpr} \cite{cakelpr} accepts the LRAT proofs of the six graphs in Theorem~\ref{thm:main} and of the second $22$-gon graph.
The $16$-gon also has a refutation of the encoding with a single unit clause, so that one case does not depend on the triangle.
\texttt{drat-trim -L} printed \texttt{s VERIFIED} for that DRAT in $6342.84$ seconds and wrote an LRAT proof of $18812690082$ bytes.
\texttt{lrat-check} printed \texttt{c VERIFIED} in $565.99$ seconds, and \texttt{cake\_lpr} printed \texttt{s VERIFIED UNSAT} in $554.40$ seconds.

Lean does not replay the SAT proof.
For $2n\in\{14,18,20,22,24\}$, and for the second $22$-gon graph, it proves that every edge has norm $1$, and that unsatisfiability of the CNF implies $\chi\ge 5$.
For the $16$-gon it proves the same two statements after dividing every file coordinate by $2\alpha$.
The two parts meet at the CNF file.
\texttt{lake exe emit\_cnf\_p14} reproduces the $14$-gon CNF byte for byte, and \texttt{lake exe emit\_cnf\_p22} does the same for the $5577$-vertex $22$-gon.
\texttt{lake exe emit\_cnf\_p22b} does the same for the $6517$-vertex $22$-gon.
\texttt{lake exe emit\_cnf\_p16} does the same for the $16$-gon.
\texttt{lake exe emit\_cnf} does the same for the $18$-gon, and \texttt{lake exe emit\_cnf\_multi} does the same for $2n\in\{20,24\}$.
\noindent The $14$-gon theorems are \texttt{P14.edge\_unit} and \texttt{P14.chi\_P14\_ge\_five}.

\noindent The $16$-gon theorems are \texttt{P16.edge\_unit} and \texttt{P16.chi\_P16\_ge\_five}.

\noindent The $18$-gon theorems are \texttt{P18.edge\_unit} and \texttt{P18.chi\_P18\_ge\_five}.

\noindent The $20$-gon theorems are \texttt{P20.edge\_unit} and \texttt{P20.chi\_P20\_ge\_five}.

\noindent The $22$-gon theorems for the graph in Theorem~\ref{thm:main} are \texttt{P22.edge\_unit} and \texttt{P22.chi\_P22\_ge\_five}.

\noindent The second $22$-gon graph has theorems \texttt{P22b.edge\_unit}\\
and \texttt{P22b.chi\_P22b\_ge\_five}.

\noindent The $24$-gon theorems are \texttt{P24.edge\_unit} and \texttt{P24.chi\_P24\_ge\_five}.

\noindent For every graph below, both \texttt{edge\_unit} and the chromatic theorem depend only on three axioms.
They are \texttt{propext}, \texttt{Classical.choice} and \texttt{Quot.sound}.

On a sumset $U+U$, a generated edge is a pair $\{p,p+u\}$ with $p,p+u\in U+U$ and $u\in U$.
An accidental edge joins two points of $U+U$ at norm $1$ whose difference does not lie in $U$.
The $4$-core is the graph left after repeatedly deleting vertices of degree less than $4$.
A deleted vertex can be coloured last, so the $4$-core is $4$-colourable if and only if the graph it came from is.
The final graphs for $2n\in\{14,16,18,20,22,24\}$ are subgraphs of these $4$-cores, and so is the second $22$-gon graph.

\begin{table}[t]
  \centering
  \small
  \begin{tabular}{rrrrrrcc}
    \toprule
    $2n$ & vertices & edges & $|U|$ & $|U+U|$ & $4$-core & Lean & $\chi\le$ \\
    \midrule
    $14$ & $225$ & $1176$ & $42$ & $841$ & $841$ & yes & $6$ \\
    $16$ & $960$ & $5248$ & $80$ & $3073$ & $2753$ & yes & $6$ \\
    $18$ & $612$ & $3042$ & $54$ & $1405$ & $1369$ & yes & $6$ \\
    $20$ & $1861$ & $11615$ & $140$ & $9501$ & $8441$ & yes & $6$ \\
    $22$ & $5577$ & $40603$ & $154$ & $11529$ & $10847$ & yes & $6$ \\
    $22$ & $6517$ & $42394$ & $154$ & $11529$ & $10055$ & yes & $6$ \\
    $24$ & $1296$ & $8616$ & $168$ & $13609$ & $13609$ & yes & $7$ \\
    \bottomrule
  \end{tabular}
  \caption{The graphs. The row with $6517$ vertices is a second graph in the $22$-gon norm, and it is not the graph named in Theorem~\ref{thm:main}. For every row the $4$-core column is the $4$-core of the unit-distance graph on $U+U$, and each final graph is a subgraph of that core. The upper bounds are Geh\'er's $6$ for $2n\le 22$ and Chilakamarri's $7$ for $2n=24$.}
  \label{tab:five}
\end{table}

\section{The graphs}\label{sec:graphs}
\subsection{$2n=14$}
Here $\alpha^3-\alpha^2-2\alpha+1=0$ and $\mathrm{sc}=1$.
The parameters are $5+\alpha-2\alpha^2$ and $-4-\alpha+2\alpha^2$.
The set $U$ consists of the vertices $\zeta^0,\dots,\zeta^{13}$ and the points $(1-\lambda)\zeta^k+\lambda\zeta^{k+1}$ for $0\le k<14$ and for these two values of $\lambda$.
Then $|U|=42$ and $|U+U|=841$.
The $4$-core of the unit-distance graph on $U+U$ is all of $U+U$, with $841$ vertices and $5908$ edges, and the graph below is a rotation-invariant subgraph of this $4$-core.
The hashed graph file does not record the parameters.
Method~(2) was run on a copy with those two parameter lines inserted after the scale line.
The graph has $225$ vertices and $1176$ edges, minimum degree $7$ and maximum degree $42$, and a triangle on vertices $0,158,186$.
Its edges use $154$ distinct unit vectors.
\begin{figure}[t]
  \centering
  \includegraphics[width=\textwidth]{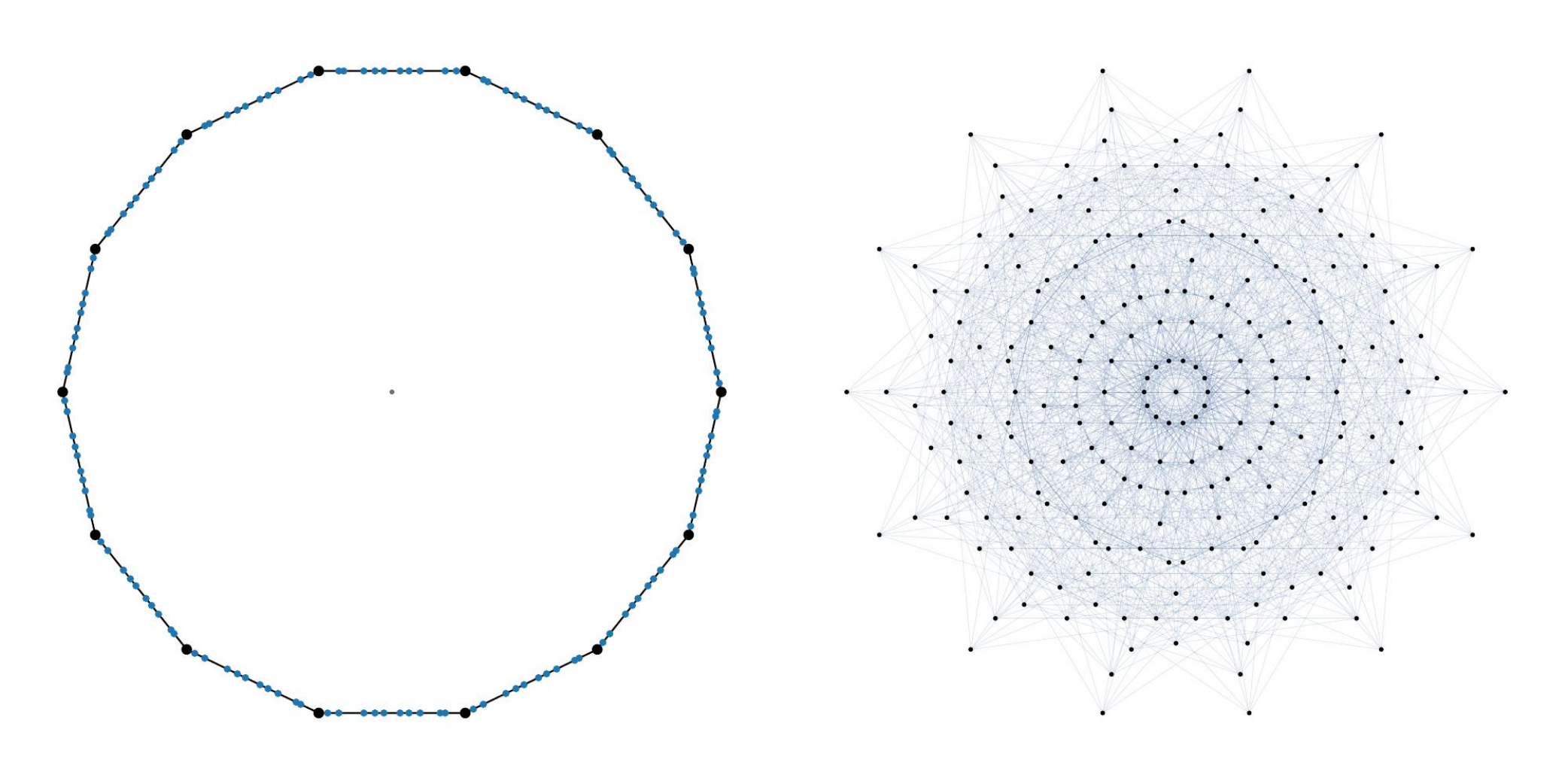}
  \caption{Left: $P_{14}$ and these $154$ unit vectors, vertices in black and the other vectors in blue. Right: the $225$ vertices and $1176$ edges.}
  \label{fig:14}
\end{figure}
Minimality was not examined.
The certificate is in Section~\ref{sec:data}.

\subsection{$2n=16$}
Let $\alpha=2\cos(\pi/8)$, so $\alpha^2=2+\sqrt2$, and let $V_k=2\alpha\zeta^k$.
The side parameters are
\[
  t=4+5\alpha-\alpha^2-\tfrac32\alpha^3,\qquad
  s=-1-2\alpha+\tfrac12\alpha^2+\tfrac12\alpha^3,
\]
together with $1-t$ and $1-s$.
The set $U$ consists of the sixteen vertices $V_k$ and the points $(1-\lambda)V_k+\lambda V_{k+1}$ for these four values of $\lambda$.
Then $|U|=80$ and $|U+U|=3073$.
The $4$-core of the $Q_{16}$-unit-distance graph on $U+U$ (which division by $2\alpha$ carries to a $P_{16}$-unit-distance graph) has $2753$ vertices, and the graph below is a rotation-invariant subgraph of this $4$-core.
It has $960$ vertices and $5248$ edges ($3968$ generated and $1280$ accidental), minimum degree $4$, maximum degree $43$, and is a union of $60$ rotation orbits.
\begin{figure}[t]
  \centering
  \includegraphics[width=\textwidth]{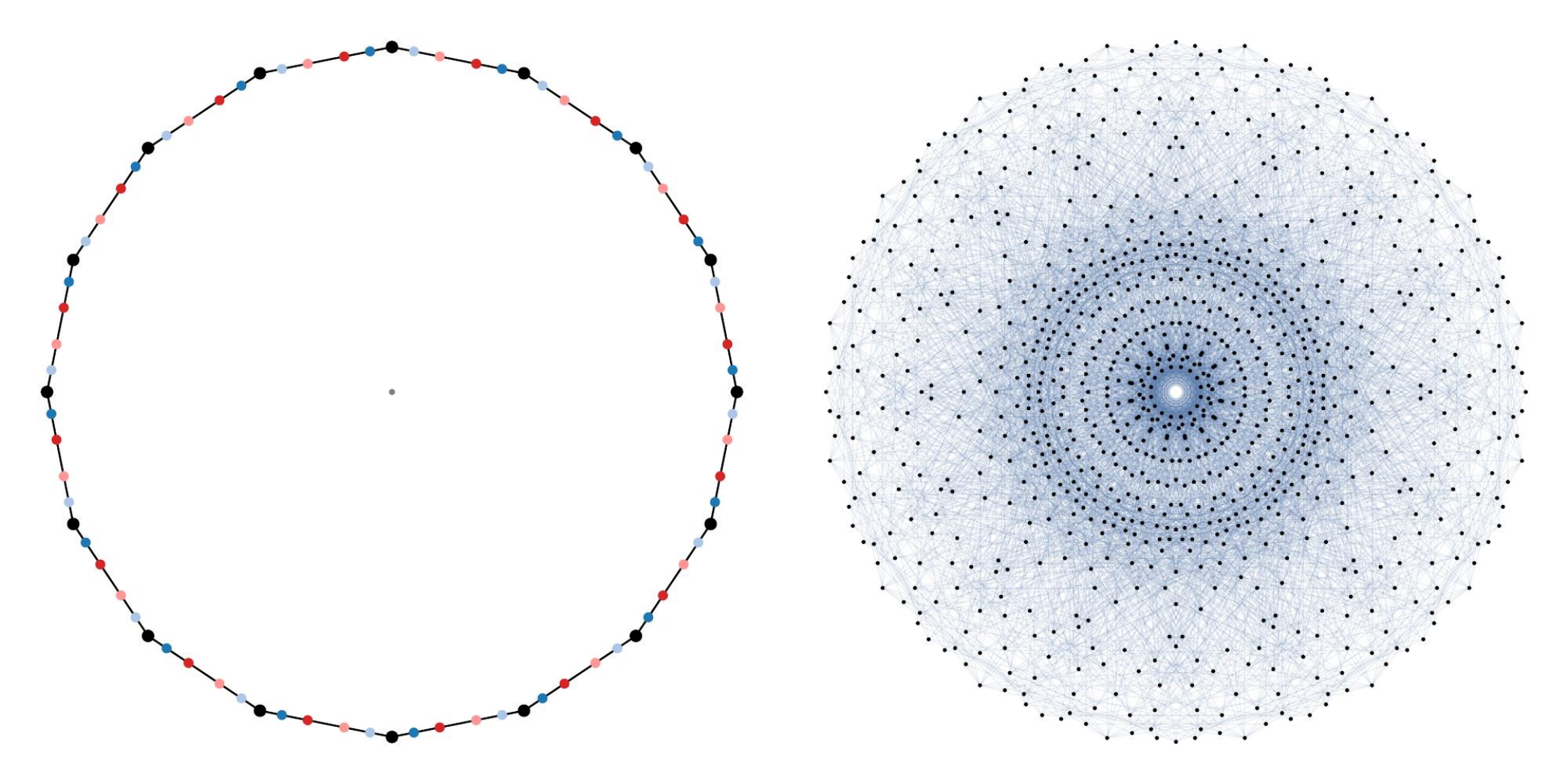}
  \caption{Left: $Q_{16}=2\alpha P_{16}$ and the $80$ points of $U$; the colours distinguish the vertices and the four side parameters. Right: the $960$ vertices and $5248$ edges.}
  \label{fig:16}
\end{figure}
It is minimal for deleting rotation orbits; vertex-minimality is not claimed.
A subgraph with $928$ vertices and $5124$ edges also has a checked DRAT and LRAT refutation.
\texttt{lake exe emit\_cnf\_p16} reproduces the triangle CNF byte for byte.
The Lean theorems are \texttt{P16.edge\_unit} and \texttt{P16.chi\_P16\_ge\_five}.
Both depend only on \texttt{propext}, \texttt{Classical.choice} and \texttt{Quot.sound}.
The formalised coordinates are the file coordinates divided by $2\alpha$, so the unit ball is $P_{16}=\operatorname{conv}\{\zeta^k\}$ rather than the similar copy $Q_{16}=\operatorname{conv}\{2\alpha\zeta^k\}$ in which the file is stored.
These Lean sources are in the deposit \cite{Pdata}.
\texttt{cake\_lpr} accepts the triangle LRAT proof.
The certificates are in Section~\ref{sec:data}.

\subsection{$2n=18$}
Parameters $\tfrac43-\tfrac13\alpha^2$ and $-\tfrac13+\tfrac13\alpha^2$, with $\alpha^3-3\alpha-1=0$ and $\mathrm{sc}=3$.
The set $U$ consists of the vertices $\zeta^0,\dots,\zeta^{17}$ and the points $(1-\lambda)\zeta^k+\lambda\zeta^{k+1}$ for $0\le k<18$ and for these two values of $\lambda$.
Then $|U|=54$ and $|U+U|=1405$.
The $4$-core of the unit-distance graph on $U+U$ has $1369$ vertices, and the graph below is a rotation-invariant subgraph of this $4$-core.
The final graph has $612$ vertices and $3042$ edges ($2430$ generated and $612$ accidental), minimum degree $4$, maximum degree $33$, in $34$ rotation orbits.
\begin{figure}[t]
  \centering
  \includegraphics[width=\textwidth]{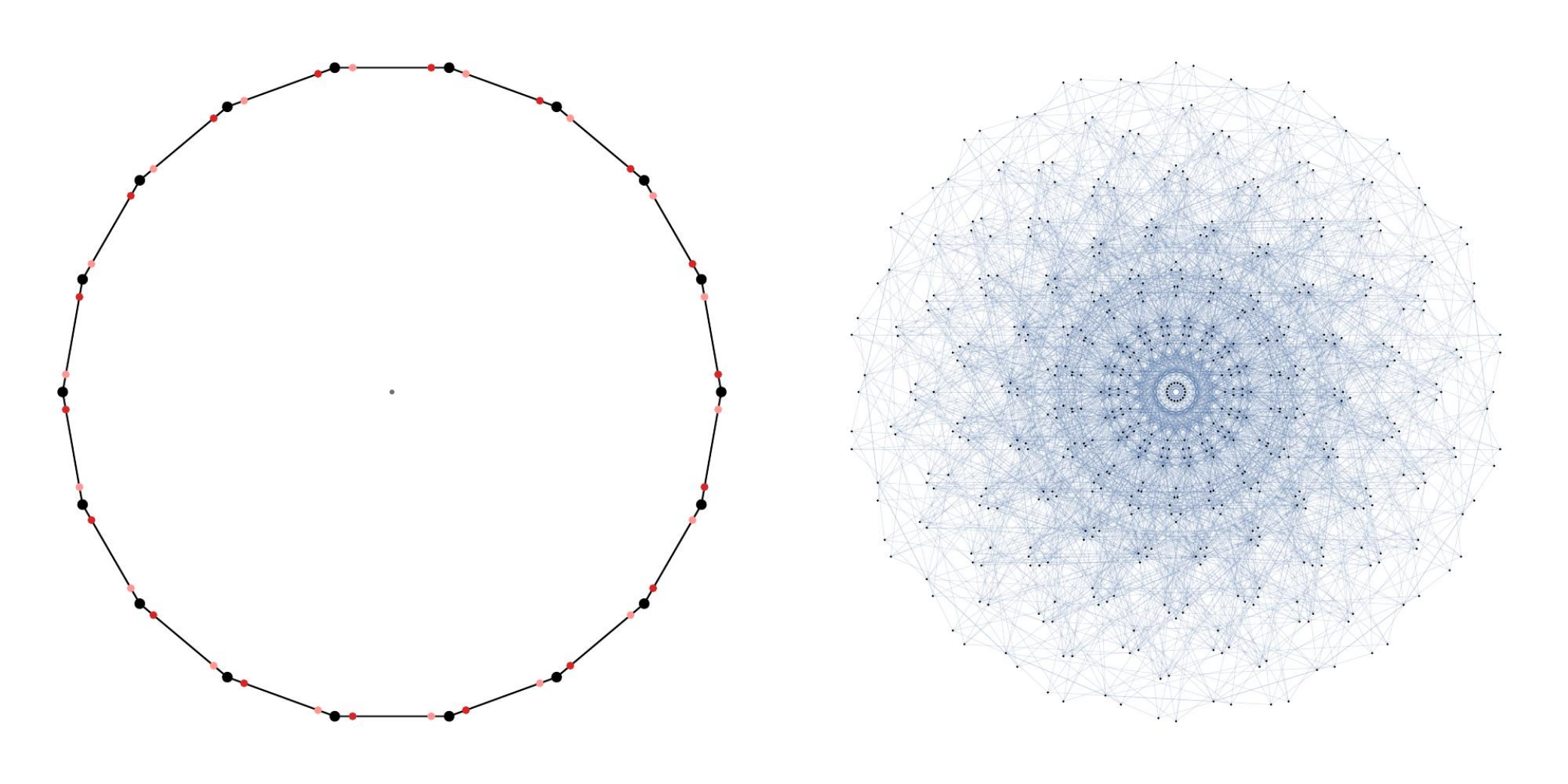}
  \caption{Left: $P_{18}$ and the $54$ points of $U$, vertices in black and one colour for each side parameter. Right: the $612$ vertices and $3042$ edges.}
  \label{fig:18}
\end{figure}
Deleting any one orbit leaves a $4$-colourable graph; vertex-minimality was not examined.
The certificate is in Section~\ref{sec:data}.

\subsection{$2n=20$}
Here $\alpha^4-5\alpha^2+5=0$ and $\mathrm{sc}=5$.
The six parameters are $8+7\alpha-2\alpha^2-2\alpha^3$, $-7-7\alpha+2\alpha^2+2\alpha^3$, $4-\alpha^2$, $-3+\alpha^2$, $-3-3\alpha+\alpha^2+\tfrac45\alpha^3$ and $4+3\alpha-\alpha^2-\tfrac45\alpha^3$.
The set $U$ consists of the vertices $\zeta^0,\dots,\zeta^{19}$ and the points $(1-\lambda)\zeta^k+\lambda\zeta^{k+1}$ for $0\le k<20$ and for these six values of $\lambda$.
Then $|U|=140$ and $|U+U|=9501$.
The $4$-core of the unit-distance graph on $U+U$ has $8441$ vertices, and the graph below is a subgraph of this $4$-core.
The final graph has $1861$ vertices and $11615$ edges ($7475$ generated and $4140$ accidental), minimum degree $4$, maximum degree $66$.
\begin{figure}[t]
  \centering
  \includegraphics[width=\textwidth]{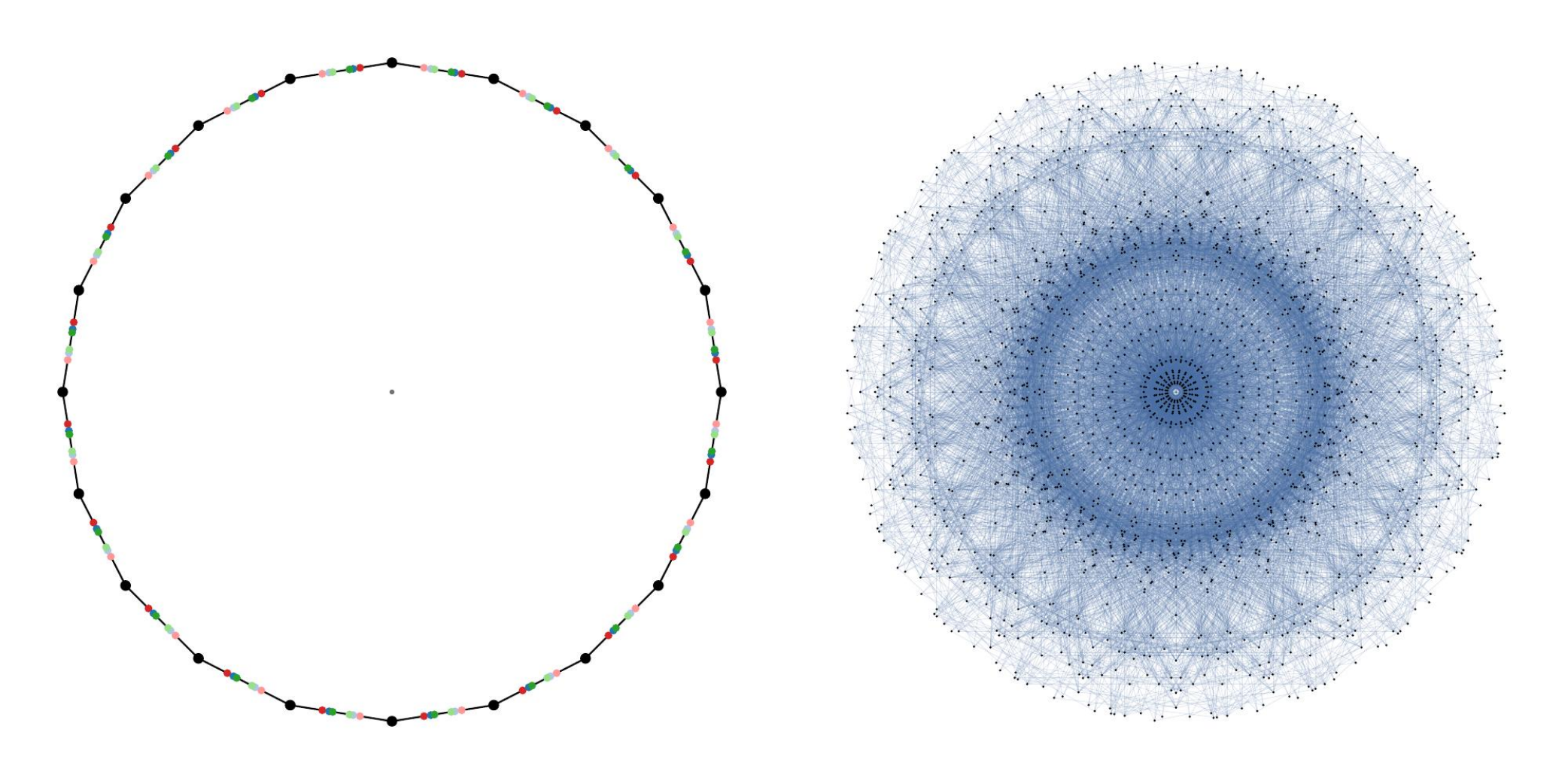}
  \caption{Left: $P_{20}$ and the $140$ points of $U$, vertices in black and one colour for each side parameter. Right: the $1861$ vertices and $11615$ edges.}
  \label{fig:20}
\end{figure}
It is not invariant under multiplication by $\zeta$. Minimality was not examined.
The certificate is in Section~\ref{sec:data}.

\subsection{$2n=22$}
Here $\alpha^5-\alpha^4-4\alpha^3+3\alpha^2+3\alpha-1=0$ and $\mathrm{sc}=1$.
The six parameters are $-2+5\alpha+11\alpha^2-\alpha^3-3\alpha^4$, $3-5\alpha-11\alpha^2+\alpha^3+3\alpha^4$, $2-4\alpha-4\alpha^2+\alpha^3+\alpha^4$, $-1+4\alpha+4\alpha^2-\alpha^3-\alpha^4$, $1+\alpha+3\alpha^2-\alpha^4$ and $-\alpha-3\alpha^2+\alpha^4$.
The set $U$ consists of the vertices $\zeta^0,\dots,\zeta^{21}$ and the points $(1-\lambda)\zeta^k+\lambda\zeta^{k+1}$ for $0\le k<22$ and for these six values of $\lambda$.
Then $|U|=154$ and $|U+U|=11529$.
The $4$-core of the unit-distance graph on $U+U$ has $10847$ vertices, and the graph below is a subgraph of this $4$-core.
The final graph has $5577$ vertices and $40603$ edges, minimum degree $4$ and maximum degree $150$.
\begin{figure}[t]
  \centering
  \includegraphics[width=\textwidth]{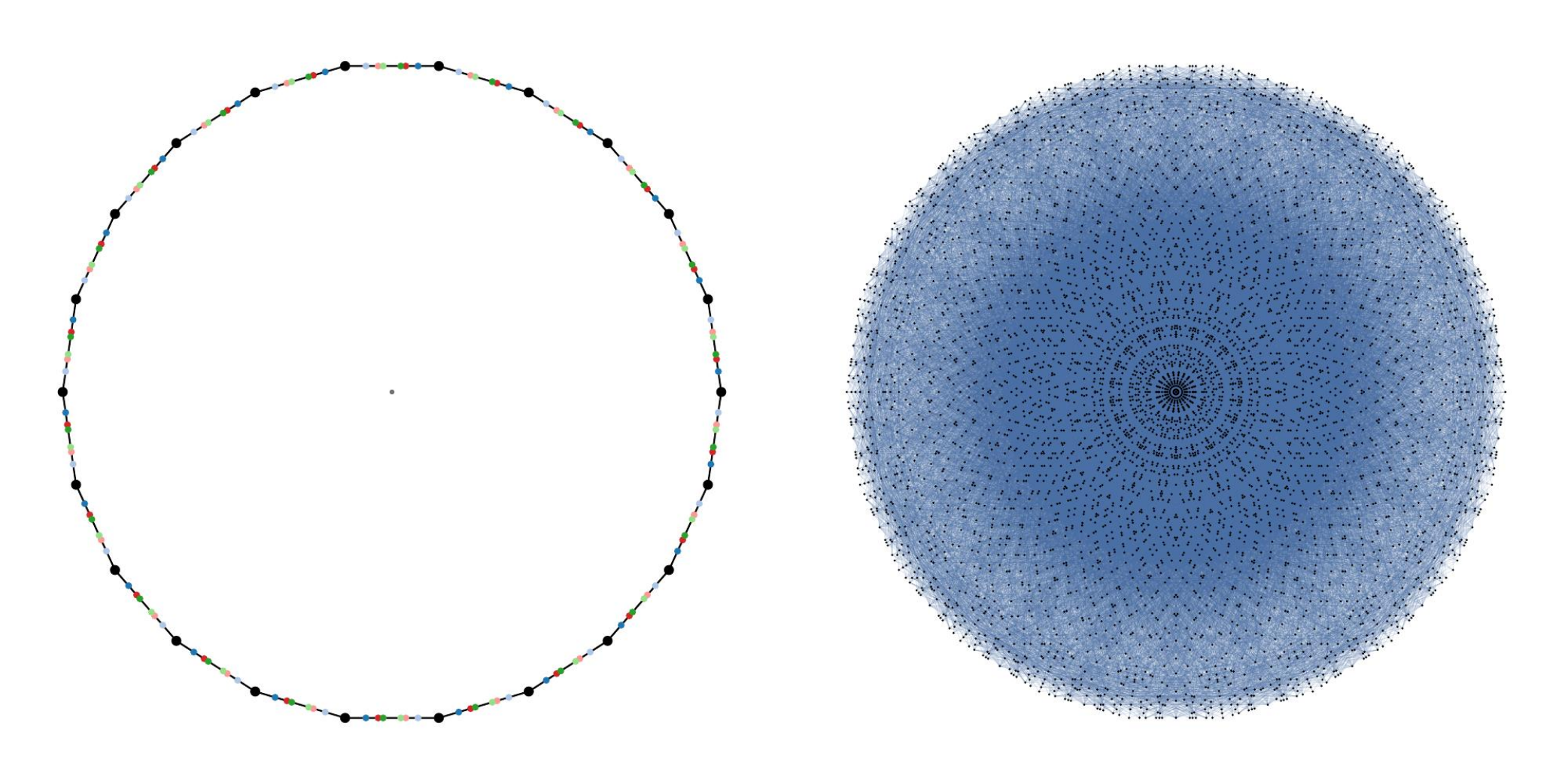}
  \caption{Left: $P_{22}$ and the $154$ points of $U$, vertices in black and one colour for each side parameter. Right: the $5577$ vertices and $40603$ edges.}
  \label{fig:22}
\end{figure}
Minimality was not examined.
The certificate is in Section~\ref{sec:data}.

\subsection{A second graph for $2n=22$}
The minimal polynomial and the scale are those of the previous subsection.
The first four parameters are unchanged.
The last two are $6+2\alpha-8\alpha^2-\alpha^3+2\alpha^4$ and $-5-2\alpha+8\alpha^2+\alpha^3-2\alpha^4$.
The set $U$ is built from the vertices and these six parameters, so again $|U|=154$ and $|U+U|=11529$.
The $4$-core of the unit-distance graph on $U+U$ has $10055$ vertices and $61292$ edges, and the graph below is a subgraph of this $4$-core.
It has $6517$ vertices and $42394$ edges, minimum degree $4$ and maximum degree $163$.
\begin{figure}[t]
  \centering
  \includegraphics[width=\textwidth]{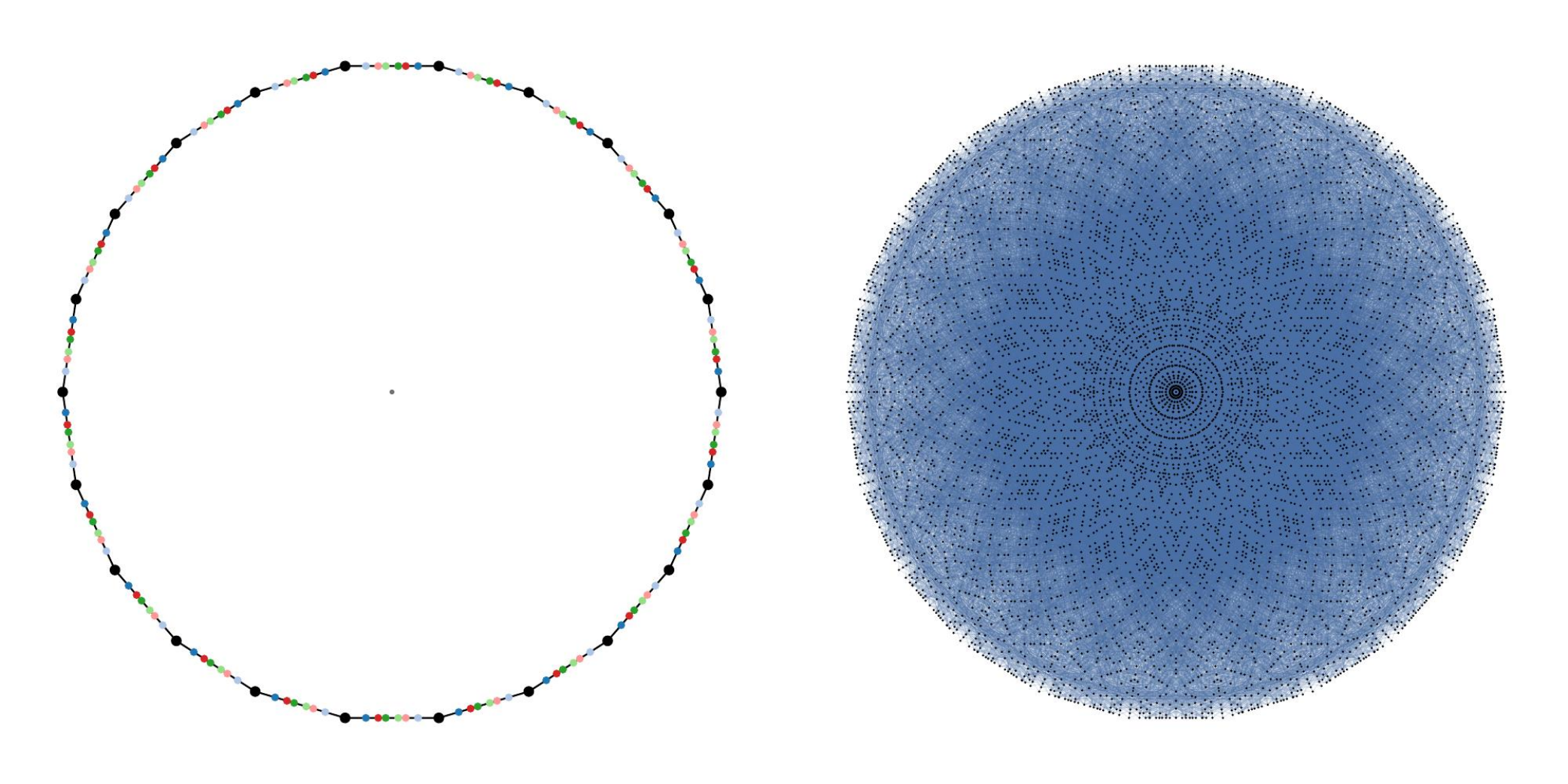}
  \caption{Left: $P_{22}$ and the $154$ points of $U$, vertices in black and one colour for each side parameter. Right: the $6517$ vertices and $42394$ edges.}
  \label{fig:22b}
\end{figure}
This is not the graph named in Theorem~\ref{thm:main}.
It is not invariant under multiplication by $\zeta$.
Minimality was not examined.
\texttt{lake exe emit\_cnf\_p22b} reproduces the triangle CNF byte for byte.
The Lean theorems are \texttt{P22b.edge\_unit} and \texttt{P22b.chi\_P22b\_ge\_five}, and both depend only on \texttt{propext}, \texttt{Classical.choice} and \texttt{Quot.sound}.
These Lean sources are in the deposit \cite{Pdata}.
The certificate is in Section~\ref{sec:data}, and the proof files are not part of that deposit.

\subsection{$2n=24$}
Here $\alpha^4-4\alpha^2+1=0$ and $\mathrm{sc}=2$.
The six parameters are $4-\alpha^2$, $-3+\alpha^2$, $1+7\alpha-2\alpha^3$, $-7\alpha+2\alpha^3$, $-\tfrac32+\tfrac12\alpha^2$ and $\tfrac52-\tfrac12\alpha^2$.
The set $U$ consists of the vertices $\zeta^0,\dots,\zeta^{23}$ and the points $(1-\lambda)\zeta^k+\lambda\zeta^{k+1}$ for $0\le k<24$ and for these six values of $\lambda$.
Then $|U|=168$ and $|U+U|=13609$.
The $4$-core of the unit-distance graph on $U+U$ has $13609$ vertices, which is all of $U+U$, and the graph below is a rotation-invariant subgraph of this $4$-core.
The final graph has $1296$ vertices and $8616$ edges ($7008$ generated and $1608$ accidental), minimum degree $6$, maximum degree $44$, in $54$ rotation orbits.
\begin{figure}[t]
  \centering
  \includegraphics[width=\textwidth]{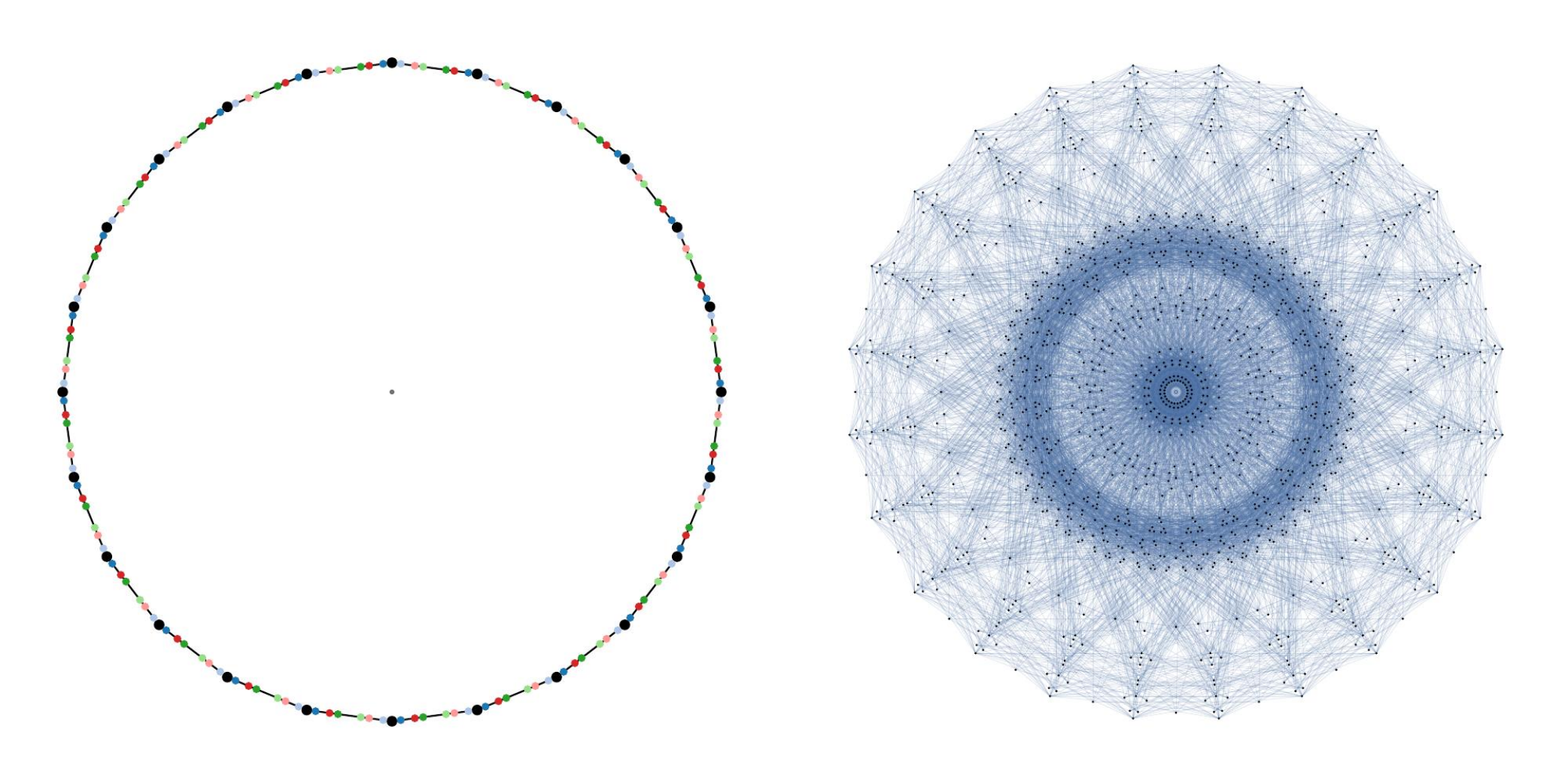}
  \caption{Left: $P_{24}$ and the $168$ points of $U$, vertices in black and one colour for each side parameter. Right: the $1296$ vertices and $8616$ edges.}
  \label{fig:24}
\end{figure}
Vertex-minimality is not claimed.
The checked proof is the core-lemma DRAT in Section~\ref{sec:data}.

\section{Data}\label{sec:data}
The $16$-gon graph and both CNF files are at doi:10.5281/zenodo.24458031 \cite{P16data}.
The triangle DRAT and LRAT proofs, and the original-encoding DRAT and LRAT proofs, are in the same record.
The graph, CNF and proof files for $2n\in\{14,18,20\}$, the $5577$-vertex $22$-gon, and the $24$-gon graph, CNF, core-lemma DRAT and LRAT, are at doi:10.5281/zenodo.24458454 \cite{Pdata}.
The Lean sources for $2n\in\{14,16,18,20,22,24\}$ and for the second $22$-gon graph are in the same record.
Each SHA-256 is printed on two lines.

\noindent\textbf{$2n=14$.}
The triangle is on vertices $0,158,186$.
CaDiCaL~3.0.1 returned UNSAT in $11.78$ seconds.
\texttt{drat-trim} printed \texttt{s VERIFIED} in $8.93$ seconds, and \texttt{lrat-check} printed \texttt{c VERIFIED} in $1.16$ seconds.
\texttt{cake\_lpr} printed \texttt{s VERIFIED UNSAT} in $0.98$ seconds.

Graph file, SHA-256\newline
\hash{daf4940e495e25ddc11d9e69c83359b9}{ce861411efb12b78040e4773b20d812d}

Triangle CNF, SHA-256\newline
\hash{f6b4a6cbe41b9317e6f911846038479c}{b8dcfc36e5f07caf9680b4218faad8c2}

DRAT file, SHA-256\newline
\hash{fdad3a3b3935d26697e3f2306ca4f67d}{7f7778d538dfb5a79d45518cc512d9d8}

LRAT file, SHA-256\newline
\hash{e45c4ef18c58be8657726013de492044}{8ad8b867e6c54bb0d7420bab85f6550a}

\noindent\textbf{$2n=16$.}
The triangle is on vertices $606,634,871$.
CaDiCaL refuted the triangle CNF in $291$ seconds. \texttt{drat-trim} accepted the DRAT proof in $281$ seconds, \texttt{lrat-check} accepted the LRAT proof, and \texttt{cake\_lpr} printed \texttt{s VERIFIED UNSAT} in $53.36$ seconds.
The original encoding, with a single unit clause, was refuted as well (CaDiCaL $5170$ seconds, \texttt{drat-trim} $5272$ seconds).
A later \texttt{drat-trim -L} printed \texttt{s VERIFIED} in $6342.84$ seconds and wrote an LRAT proof.
\texttt{lrat-check} printed \texttt{c VERIFIED} in $565.99$ seconds, and \texttt{cake\_lpr} printed \texttt{s VERIFIED UNSAT} in $554.40$ seconds.

Graph file, SHA-256\newline
\hash{97ae6cf2a4b8cd5039f78f1a664d6fbe}{79369f8f63dcf2624a7dabdb6e4616f8}

Triangle CNF, SHA-256\newline
\hash{602871519fe0dce62d71ba44dcd9c3d8}{89740226a99e1f932072a56bf7752d24}

Original CNF, SHA-256\newline
\hash{2cf574dd4d2c6ca0dddcba5e47d28f5d}{30d4f7167e11513cb41fba2b759042db}

Triangle DRAT, SHA-256\newline
\hash{e80589e2b809e6772e856130c572ce69}{80f9b094f275829b2972ef09af7e96df}

Original DRAT, $16807562552$ bytes, SHA-256\newline
\hash{92dbfc1dc255741222202358e9dc4ce1}{4e45ce2adc073501e89c0eaaacbd066f}

Original LRAT, $18812690082$ bytes, SHA-256\newline
\hash{36e4d14c9a2956e7cdc8b2b04006654e}{9864d52c3c303113ec1954791c73daaf}

Triangle LRAT, $1857281172$ bytes, SHA-256\newline
\hash{3c934b5f23d1264cde94dc4ab469ae49}{1dc6c8d97876abd6c64c41dfab109515}

\noindent\textbf{$2n=18$.}
The triangle is on vertices $525,582,610$.
CaDiCaL refuted the triangle CNF in $293.35$ seconds.
\texttt{drat-trim} took $269.44$ seconds, \texttt{lrat-check} $18.82$ seconds, and \texttt{cake\_lpr} $53$ seconds.
The interval check deviates from $1$ by at most $6.000\times 10^{-58}$, against separation at least $8.333\times 10^{-6}$.

Graph file, SHA-256\newline
\hash{e1ff14fe72aae17b4a59a65ff3df8731}{c797ad512baa00491f9808f05f0b41bd}

Triangle CNF, SHA-256\newline
\hash{1ee31b5aae0e5f477f438445100357c5}{ed28c630839d8fb82a6d7ce29611a8d3}

DRAT file, $1355604349$ bytes, SHA-256\newline
\hash{466e536e36e92ad5013266c30ec1bed5}{91f7928feb4a706b31f18f2e95f675ae}

LRAT file, $1503274894$ bytes, SHA-256\newline
\hash{815b3e55ccde0e1fa03b98fb6913058b}{68a0d8c228583965069f0052ee24e5d4}

\noindent\textbf{$2n=20$.}
The triangle is on vertices $547,568,582$.
CaDiCaL took $2990.16$ seconds.
\texttt{drat-trim} took $2703.51$ seconds, \texttt{lrat-check} $118.2$ seconds, and \texttt{cake\_lpr} $253$ seconds.
The interval deviation from $1$ is at most $1.277\times 10^{-57}$, against separation at least $7.337\times 10^{-11}$.

Graph file, SHA-256\newline
\hash{e20e6618b4c401a0da7c09f2ec461fca}{18fab44b7e34032cc4e471e4f82f9c06}

Triangle CNF, SHA-256\newline
\hash{cd1c0cbefa7ad40d456edfbcb1eabaf0}{3f73e1d918e9ac20147a3e990fb4161d}

DRAT file, $6002775922$ bytes, SHA-256\newline
\hash{b94be8e6564dac30dfef8e67373b905c}{f67dcbae96710cd3ad9c939c91b480ad}

LRAT file, $6797768876$ bytes, SHA-256\newline
\hash{a6f617527bbffd7ae1b636313507972b}{31a3a80ee2fa1c4d2241c0c8dc39c15a}

\noindent\textbf{$2n=22$.}
The triangle is on vertices $2222,2225,2795$.
CaDiCaL returned UNSAT in $613.25$ seconds.
\texttt{drat-trim -l} printed \texttt{s VERIFIED} in $381.5$ seconds and wrote a core-lemma DRAT of $812850208$ bytes.
A second \texttt{drat-trim} printed \texttt{s VERIFIED} in $253.3$ seconds and wrote the LRAT proof.
\texttt{lrat-check} printed \texttt{c VERIFIED} in $25.55$ seconds, and \texttt{cake\_lpr} printed \texttt{s VERIFIED UNSAT} in $58.71$ seconds.
The interval deviation from $1$ is at most $1.858\times 10^{-57}$, against separation at least $4.567\times 10^{-11}$.

Graph file, SHA-256\newline
\hash{efe67e0f7d7f15d12bbabb2571e01c24}{3c4b3637c70784fc7e495f04328e21bc}

Triangle CNF, SHA-256\newline
\hash{7f2f69f57dcdbb71fcf65b25dc0ce328}{573a926aaeb13edff5b65febee5f0765}

Original DRAT, $2411499188$ bytes, SHA-256\newline
\hash{8f01522b3400d6d5c45e190d70d57f69}{23a83f070e1d0c879151c39ef9304b69}

Core-lemma DRAT, SHA-256\newline
\hash{18166372f15136ccf4a1b2b214dca8e2}{2cb03486f3049c0f77246091b29d3b04}

LRAT file, $2305010264$ bytes, SHA-256\newline
\hash{76888b0393960388a0efcf8f90eefba8}{640e209fabe71f99b5a8f2ec82af0d79}

\noindent\textbf{A second graph for $2n=22$.}
The triangle is on vertices $732,734,3253$.
CaDiCaL returned UNSAT in $11650.84$ seconds.
\texttt{drat-trim -l} printed \texttt{s VERIFIED} in $13283.28$ seconds.
It wrote a core-lemma DRAT of $10429730245$ bytes.
A second \texttt{drat-trim} printed \texttt{s VERIFIED} in $8262.86$ seconds and wrote the LRAT proof.
\texttt{lrat-check} printed \texttt{c VERIFIED} in $1574.98$ seconds, and \texttt{cake\_lpr} printed \texttt{s VERIFIED UNSAT} in $1736.22$ seconds.
The interval deviation from $1$ is at most $1.858\times 10^{-57}$, against separation at least $4.567\times 10^{-11}$.
These files are not in the deposit \cite{Pdata}.

Graph file, SHA-256\newline
\hash{6d37fbd860b64d47086124277c1306a3}{7379a9d52b0352b709d03f9118946a3e}

Triangle CNF, SHA-256\newline
\hash{b22275e6a030f79bfaff8e1d5f263428}{4a52f797dd40391a55e046e6ed3dd221}

Original DRAT, $36018514927$ bytes, SHA-256\newline
\hash{94da59170b4d68490fa16fc2f624b1da}{8517beceb0e600189d8316324b3cf692}

Core-lemma DRAT, SHA-256\newline
\hash{e76984165eee6e99d5de892a0bbcc4a0}{077e9b30a56621b1a958bfef0cad3a51}

LRAT file, $54165026725$ bytes, SHA-256\newline
\hash{debd38e79a06c4eb07fd7b5545033cd5}{36a27fc72fb5c3d7a3f533f8d4381aa0}

\noindent\textbf{$2n=24$.}
The triangle is on vertices $869,884,902$.
CaDiCaL took $3759.28$ seconds and wrote a DRAT proof of $17377116364$ bytes.
That file was deleted after \texttt{drat-trim -l} produced a core-lemma proof of $5472492911$ bytes, which \texttt{drat-trim} verified in $1688.04$ seconds.
\texttt{lrat-check} took $180.48$ seconds and \texttt{cake\_lpr} $665$ seconds.
The interval deviation from $1$ is at most $1.250\times 10^{-57}$, against separation at least $1.240\times 10^{-8}$.
The original DRAT was read twice before deletion. Both reads gave SHA-256\newline
\hash{14a4bdad14fd798d35d06ae1fa24616f}{bf40b3c17d87f20afd7d011522581658}

Graph file, SHA-256\newline
\hash{29a805648abcbd9e037f9db284597591}{26ba60c4a8f33bd59d0bf8dec87f8881}

Triangle CNF, SHA-256\newline
\hash{37af375f5d10780f41e8f982308d5489}{aa867f786527a3ecfdee739971e8337f}

Core-lemma DRAT, SHA-256\newline
\hash{ec8465bcaeab704429080b91eb8db642}{48daebb5dec2c7ebb5d07b1648877b01}

LRAT file, $16878314761$ bytes, SHA-256\newline
\hash{ffefe38ca3487291b0cc52af4de2d9b0}{6bfbd7d31c81c7b6ecd89ff4ca5ea165}

\section*{Acknowledgements}
I thank Panna Geh\'er for correspondence on the status of the problem for the regular sixteen-gon, and for pointing out the preprint of OpenAI on the Euclidean plane \cite{OpenAI}.
That preprint concerns the Euclidean norm only.

\section*{Declaration of AI use}
AI assistants were used under the author's direction in writing the search and verification programs, the Lean~4 formalisation, and the manuscript text.
The author verified every mathematical argument and written proof line by line, independently ran all exact-arithmetic, SAT, DRAT, LRAT, \texttt{cake\_lpr} and Lean~4 checks reported in Sections~\ref{sec:chain} and~\ref{sec:data}, and takes full responsibility for the content of the paper.

\end{document}